\documentstyle[11pt]{article}
\def\pl{\partial}

\def\*{\raisebox{.5mm}{*}}

\def\div{{\,\rm div\,}}

\def\pa{\partial}

\def\R{I\!\!R}
\def\A{{\cal A}}
\def\H{{\cal H}}

\def\Ga{\Gamma}

\def\<{\langle}
\def\>{\rangle}

\def\Om{\Omega}

\def\det{\mbox{det\,}}

\def\bma{\left[\begin{array}}
\def\ema{\end{array}\right]}
\def\bda{\left|\begin{array}}
\def\eda{\end{array}\right|}
\def\be{\begin{equation}}
\def\ee{\end{equation}}

\newtheorem{thm}{{}\hskip\parindent Theorem}[section]
\newtheorem{lem}{{}\hskip\parindent Lemma}[section]
\newtheorem{pro}{{}\hskip\parindent Proposition}[section]
\newtheorem{exl}{{}\hskip\parindent Example}[section]

\newtheorem{dfn}{{}\hskip\parindent Definition}[section]
\newtheorem{rem}{{}\hskip\parindent Remark}[section]

\large
\title{Special uniform decay rate of  local energy for the wave equation with variable coefficients on an exterior domain}
\date{}
\author{Zhen-Hu Ning\thanks{Corresponding author, E-mail address: nzh41034@163.com}, Fengyan Yang and Xiaopeng Zhao}
\begin{document}
\maketitle

\footnote{Zhen-Hu Ning,
 Faculty of Information Technology, Beijing University of Technology, Beijing, 100124, China. E-mail address: nzh41034@163.com.

\ \ Fengyan Yang,
Key Laboratory of Systems and Control, Institute of Systems Science, Academy of Mathematics and
Systems Science, Chinese Academy of Sciences, Beijing, 100190, China.
E-mail address: yangfengyan12@mails.ucas.edu.cn.

\ \
Xiaopeng Zhao,
School of Science, Jiangnan University,
Wuxi, Jiangsu, 214122, China.
E-mail address: zhaoxiaopeng@jiangnan.edu.cn.
}
\footnotetext{This work is supported by the National Science Foundation of China, grants  no.61473126 and no.61573342, and Key Research Program of Frontier Sciences, CAS, no. QYZDJ-SSW-SYS011.}
\begin{quote}
\begin{small}
{\bf Abstract} \,\, We consider the wave equation with variable
coefficients on an exterior domain in $\R^n$($n\ge 2$).
  We are interested in finding
a special uniform decay rate of local energy  different from  the constant coefficient wave equation.
 More concretely, if the  dimensional $n$ is even, whether the uniform decay rate of local energy  for
 the wave equation with variable
coefficients can break through the limit of polynomial and reach exponential; if the  dimensional $n$ is odd,  whether the uniform decay rate of local energy  for
 the wave equation with variable
coefficients can  hold exponential as the constant coefficient wave equation .

\quad \ \ We propose a cone and establish Morawetz's multipliers in a version of the Riemannian geometry to derive uniform  decay of local energy for the wave equation with variable
coefficients. We find that the  cone with  polynomial growth is closely related to  the uniform decay rate of the local energy.  More concretely, for radial solutions, when the cone has polynomial of degree $\frac{n}{2k-1}$ growth, the uniform  decay rate of local energy  is exponential; when the cone has polynomial of degree $\frac{n}{2k}$ growth, the uniform  decay rate of local energy  is polynomial at most. In addition, for general solutions, when the cone has polynomial of degree $n$ growth, we prove that the uniform  decay rate of local energy  is exponential under suitable Riemannian metric.  It is worth pointing out that such results are independent of the parity of the dimension $n$, which is the main difference with the constant coefficient wave equation. Finally, for general solutions, when the cone has polynomial of degree $m$ growth, where $m$ is
any positive constant, we prove that the uniform decay rate of the local energy is of primary polynomial under suitable Riemannian metric.
\\[3mm]
\\[3mm]
\end{small}
\end{quote}
\vskip .5cm
\tableofcontents
\def\theequation{1.\arabic{equation}}
\setcounter{equation}{0}
\section{Introduction}
\vskip .2cm
 \quad  \ \ Let $O$ be the original point of $\R^n$ ($n\ge 2$) and
 \be r(x)= |x|,\quad  x\in\R^n\ee
be  the standard distance function of $\R^n$.
Moreover, let $\<\cdot,\cdot\>$, $\div$, $\nabla$, $\Delta$ and $I_n=(\delta_{i,j})_{n\times n}$ be the standard inner product of $\R^n$,  the standard divergence operator of $\R^n$, the standard gradient operator of $\R^n$, the standard Laplace operator of $\R^n$ and the unit matrix.

Let $\Omega$ be an exterior domain in $\R^n$ with a compact smooth boundary
$\Ga$. We assume the original point $O\notin \overline{\Omega}$.
Let
\be \label{cpde.3.1}r_0= \inf_{x\in \Ga}|x|,\quad r_1= \sup_{x\in \Ga}|x|.\ee
Then $r_1\ge r_0>0$.
$$ $$

 We consider the following system.  
\begin{equation}
\label{wg.1} \cases{u_{tt}-\div A(x) \nabla u=0\qquad (x,t)\in \Omega\times
(0,+\infty),\cr
 u\large|_{\Ga}=0\qquad t\in(0,+\infty),
\cr u(x,0)=u_0(x),\quad u_t(x,0)=u_1(x)\qquad x\in \Omega,}
\end{equation}
where $A(x)=(a_{ij}(x))_{n\times n}$  is a symmetric, positive definite matrix  for each $x\in\R^n$ and $a_{ij}(x)$ ($1\leq i,j \leq n$) are smooth functions on $\R^n$.

For $a>r_1$, the local energy for the system (\ref{wg.1}) is defined by

\begin{equation}
\label{wg.2}  E(t,a)=\frac12\int_{\Omega(a)}\left(u_t^2+\sum_{i,j=1}^na_{ij}u_{x_i}u_{x_j}\right)dx,
\end{equation}
where $\Om(a)=\{x|x\in \Om,|x|\leq a\}$. In this paper, we are interested in the uniform decay rate of the local energy $E(t,a)$.

 If $A(x) \equiv I_n$, the system (\ref{wg.1}) is known as constant coefficient. 
In the case of constant coefficient,  this problem has a long history and a wealth of results were obtained, see for example  
\cite{w25,36,w11,w4,w16,w6,w7,3,6,w5,5} and the reference cited therein. 
We recall that Wilcox\cite{w5} established the uniform decay of local energy  with a spherical obstacle by analyzing the explicit expression for the solution obtained by separation of variables. In \cite{w6}, by using the multiplier method, Morawetz proved that if the obstacle is star-shaped, the rate of decay is $1/\sqrt{t}$. Afterwards, this famous result has been constantly improved by \cite{w25,w11,w16,w7,5} and many others.

The system (\ref{wg.1}) is referred to as variable coefficients, where $A(x)$ is given by
the material in application. For $A(x)$ which is assumed to be constant near infinity,  some decay
estimates of the local energy were presented in \cite{w1} and \cite{2}. To be specific, By proving
that there are no resonances in a region close to the real axis, Burq\cite{w1} obtained the logarithmic decay estimate of the local energy for smooth data without non-trapping assumption. While Vainberg\cite{2} studied the uniform decay of system (\ref{wg.1}) for nonsmooth data under the non-trapping assumption. 
As it can help us understand the problem from the physical point of view, the non-trapping assumption has already been widely used for the dispersive estimates(see \cite{w20,1ww,ww1,w31,w2,w3,76,w32} and so on). It says: $\lim_{t\rightarrow +\infty} \rho(\gamma (t))=+\infty,$ for any $x\in \R^n$ and any geodesic $\gamma (t)$ starting at $x$, where $\rho(x)$ is the geodesic distance function. Since the geodesic   depends on the nonlinear ODE, the non-trapping assumption is hard to check. In fact, it just
changs the dispersive problem of linear PDEs with variable coefficients into
some problems of nonlinear ODEs. And some useful criterions are further needed to make the the non-trapping assumption checkable.

As is known, the multiplier method is a simple and effective tool to deal with the energy  estimate on PDEs. In particular,  the celebrated Morawetz's multipliers introduced by \cite{w6} have been extendedly used for studying the energy decay of the wave equation with
constant coefficients, see \cite{w25,w11,w16,w7,5} and many others. Therefore, one purpose
of this paper is to establish Morawetz＊s multipliers in a version of the Riemannian geometry.



Define \be\label{g} g=A^{-1}(x)\quad\mbox{for}\quad x\in\R^n\ee
as a Riemannian metric on $\R^n$ and consider the couple
$(\R^n,g)$ as a Riemannian manifold with an inner product
\be\<X,Y\>_g=\<A^{-1}(x)X,Y\>,\quad X,Y\in\R^n_x,x\in\R^n.\ee


Let $D$ denote the Levi-Civita connection of the metric $g$ and H be a vector field,
then the covariant differential $DH$ of the vector field H is a tensor field of rank 2 as follow:
\be DH(X,Y)(x)=\<D_YH,X\>_g(x)\quad  X,Y\in\R^n_x, x\in\R^n.\ee

Let $\nabla_g$ be the gradient operator of the Riemannian manifold $(\R^n,g)$, then
\be
\nabla_gu=A(x)\nabla u,\quad x\in\R^n.\ee

Denote
\be|X|_g^2=\<X,X\>_g,\quad X\in\R^n_x,x\in\R^n.\ee
Then \be|\nabla_gu|^2_g= \<\nabla u,A(x)\nabla u\>=\sum_{i,j=1}^na_{ij}u_{x_i}u_{x_j},\quad x\in\R^n.\ee

Hence, the formula (\ref{wg.2}) can be rewritten as
 \be E(t,a)=\frac12\int_{\Omega(a)}\left(u_t^2+|\nabla_gu|^2_g\right)dx.\ee

 In this paper, we propose a cone structure,
  which connects the Riemannian metric $\<.,.\>_g$ and the standard dot metric $\<.,.\>$,
  to study the uniform decay of local energy for the system (\ref{wg.1}).  We are interested in finding
a special decay rate of local energy  different from  the constant coefficient wave equation.
 More concretely, for the even dimensional space, whether the uniform decay rate of local energy  for
 the wave equation with variable
coefficients can break the limit of polynomial and reach exponential; for the odd dimensional space, whether the uniform decay rate of local energy  for
 the wave equation with variable
coefficients can  hold exponential as the constant coefficient wave equation .

The organization of our paper goes as follows. In Section 2, we will state our main results. Then some multiplier identities and key lammas for problem (\ref{wg.1}) will be present in Section 3. We will show the well-posedness and propagation property of system (\ref{wg.1}) in Section 4. And we will discuss the local energy decay for radial solutions of system (\ref{wg.1}) in Section 5. The technical details of the proofs of the decay results for general solutions will be given in the last two sections.

 \vskip .5cm
\def\theequation{2.\arabic{equation}}
\setcounter{equation}{0}
\section{Main results}
\vskip .2cm
\subsection{Well-posedness  and   propagation property}

\quad \ \
 Let $r_0$ be given by (\ref{cpde.3.1}), we define \be \label{wg.4_1} F(y)= \sup\left\{\sqrt{\left\<\frac{\partial }{\partial r},A(x) \frac{\partial }{\partial r}\right\>}(x) \quad \Big|\quad x\in\R^n, |x|=y\right\}, \quad y\ge r_0.\ee
{\bf Assumption (A)}\,\,\,$F(y)$ satisfies
\begin{equation}
\label{wg.4}\int_{r_0}^{+\infty}\frac{1}{F(y)}dy=+\infty.
\end{equation}

Denote
\be C^{\infty}_1(\Om)=\{w\in C^{\infty} (\overline{\Om}),\  u|_{\Ga}=0\quad and \quad \int_{\Omega}\sum_{i,j=1}^na_{ij}w_{x_i}w_{x_j}dx<+\infty\}.\ee
\be C^{\infty}_2(\Om)=\{w\in C^{\infty} (\overline{\Om}),\  u|_{\Ga}=0\quad and \quad \int_{\Om} \left( \sum_{i,j=1}^na_{ij}w_{x_i}w_{x_j}+(\div A(x) \nabla w)^2\right)  dx <+\infty \}.\ee
Let $\widehat{H}^1_{0}(\Om)$ be the closure of $C^{\infty}_1(\Om)$ with respect to the tolopogy
\be  \|w(x)\|^2_{\widehat{H}^1_{0}(\Om)}=\int_{\Omega}\sum_{i,j=1}^na_{ij}w_{x_i}w_{x_j}dx.\ee
and  $\widehat{H}^2_{0}(\Om)$ be the closure of $C^{\infty}_2(\Om)$ with respect to the tolopogy
\be  \|w(x)\|^2_{\widehat{H}^2_{0}(\Om)}=\int_{\Om} \left( \sum_{i,j=1}^na_{ij}w_{x_i}w_{x_j}+(\div A(x) \nabla w)^2\right)  dx.\ee
The well-posedness of the system (\ref{wg.1}) is derived as follows.
\begin{thm}\label{wg.5} Let Assumption ${\bf (A)}$ hold true. Then, for any initial datum $(u_0,u_1)\in \widehat{H}^1_{0}(\Om)\times L^2(\Om)$ , there exists a unique solution $u$ of the system (\ref{wg.1}) satisfying $u_t\in C([0,+\infty),L^2(\Om))$ and $u\in C([0,+\infty),\widehat{H}^1_{0}(\Om))$.

 Moreover, if $(u_0,u_1)\in \widehat{H}^2_{0}(\Om)\times  (\widehat{H}^1_{0}(\Om)\cap  L^2(\Om) )$, then
the unique solution $u$ satisfies
 $u_t\in C([0,+\infty),(\widehat{H}^1_{0}(\Om)\cap  L^2(\Om) ))$ and $u\in C([0,+\infty), \widehat{H}^2_{0}(\Om))$.\end{thm}

Define the energy of the system (\ref{wg.1}) by
\begin{equation}
\label{wg.3}  E(t)=\frac12\int_{\Omega}\left(u_t^2+|\nabla_gu|^2_g\right)dx.
\end{equation}
From Theorem \ref{wg.5}, if Assumption ${\bf (A)}$ holds true, we have $E(t)=E(0), \quad \forall  t>0$.\\

The finite speed of propagation property of the wave with variable
coefficients can be stated as follows.

\begin{thm}\label{wg.c5} Let Assumption ${\bf (A)}$ hold true and let initial datum $u_0,u_1$ satisfy
\be   u_0(x)=u_1(x)=0, \quad |x|\ge R_0,\ee
where $R_0>r_1$ is a  constant.
Then
\be u(x,t)=0,\quad \textmd{for}~|x|\ge z +R_0,\ee where $z$ satisfies \be\int_{R_0}^{z+R_0}\frac{1}{F(y)}dy=t.\ee\end{thm}


\vskip .2cm
\subsection{Main assumptions}
\begin{dfn}\,\,\, We say $g$ is a {\bf cone near infinity} if there exists $M\ge r_0$, where $r_0$ is given by (\ref{cpde.3.1}), such that
\begin{equation}
\label{pde.1}A(x)x=\frac{1}{\phi^2(r)}x,\quad \textmd{for} \ \  |x|\ge M,
\end{equation}
where $\phi(y)\in C^\infty([M,+\infty))$  satisfies
\be \phi(y)>0,\quad y\in[M,+\infty),\ee
\be\int_{M}^{+\infty}\phi(y)dy=+\infty.\ee
If $M=r_0$, $g$ is called a {\bf cone}.
\end{dfn}

\begin {rem}
If  $g$ is a  cone near infinity, with (\ref{wg.4_1}) we have
\be   F(r)=\sqrt{\left\<\frac{\partial }{\partial r},A(x) \frac{\partial }{\partial r}\right\>}= \frac{1}{\phi(r)},\quad |x|\ge M,  \ee
 where $M\ge r_0$ is a positive constant. Then
     \be \int_{r_0}^{+\infty}\frac{1}{F(y)}dy\ge\int_{M}^{+\infty}\frac{1}{F(y)}dy =\int_{M}^{+\infty}\phi(r)dy=+\infty, \ee
Assumption ${\bf (A)}$ holds true. From Theorem \ref{wg.5}, the system (\ref{wg.1}) is well-posed.
\end {rem}

Next, we give some examples of the cone $g$.
\begin {exl} Let  $M\ge r_0$ be a positive constant and let
\be A(x)=\frac{1}{\phi^2(r)}I_n, \quad \textmd{for} \ \ |x|\ge M.\ee Then
  \be A(x)x=\frac{1}{\phi^2(r)}x,  \quad \textmd{for} \ \ |x|\ge M.\ee
Thus, $g$ is a cone near infinity. And if $M=r_0$, $g$ is a cone.
  \end {exl}

\begin {exl} Let  $M\ge r_0$ be a positive constant and let
 \be A(x)=\frac{x\otimes x}{\phi^2(r)|x|^2}+\left(Q(x)-\frac{ Q(x) x\otimes x}{|x|^2}\right),\quad \textmd{for} \ \ |x|\ge M,\ee
where $Q(x)=(q_{i,j})_{n\times n}(x)$ is a symmetric,  definite matrix  for each $|x|\ge r_0$ and $q_{i,j}(x)(1\leq i,j \leq n)$ are smooth functions.
Then
  \be A(x)x=\frac{1}{\phi^2(r)}x+ (Q(x)x- Q(x)x)=\frac{1}{\phi^2(r)}x, \quad \textmd{for} \ \ |x|\ge M.\ee
   Thus, $g$ is a cone near infinity. And if $M=r_0$, $g$ is a cone.
\end {exl}

Let $(r,\theta)=(r,\theta_1,\theta_2,\cdots,\theta_{n-1})$ be the polar coordinates of $x\in \R^n$ in the Euclidean metric.
Then, for the cone $g$  and $|x|\ge r_0$,
\be\left\<\frac{\partial}{\partial\theta_i},\frac{\partial}{\partial r}\right\>_g=\left\<\frac{\partial}{\partial\theta_i},A^{-1}(r,\theta)\frac{\partial}{\partial r}\right\>=0,\quad 1\le i\le n-1, \ee
\be\left\<\frac{\partial }{\partial r},\frac{\partial }{\partial r}\right\>_g=\left\<\frac{\partial }{\partial r},A^{-1}(r,\theta)\frac{\partial }{\partial r}\right\>=\phi^2(r).\ee
Then the cone $g$ in the coordinates $(r,\theta)$ has the following form:
\be\label{pde.1_1} g= \phi^2(r)dr^2+\sum_{i,j=1}^{n-1}\gamma_{ij}(r,\theta)d\theta_id\theta_j,\quad \gamma_{ij}= \<A^{-1}(r,\theta)\frac{\partial}{\partial\theta_i},\frac{\partial}{\partial\theta_j}\>,\quad |x|\ge r_0. \ee

We define
\be \label{cpde.2.1} \rho(x)= \int_{r_0}^{r(x)}\phi(y)dy+c(r_0),\quad|x|\ge r_0,\ee
where $c(r_0)$ is any positive constant. Specially, if
\be  \phi(r)=m r^{m-1}, \quad  |x|\ge r_0,  \ee
where  $m>0$ is a constant. We set
 $c(r_0)=r^{m}_0$,
  then \be\label{cpde.2.2} \rho(x)=r^m,\quad  |x|\ge r_0. \ee

 Noting that $d\rho= \phi(r)dr$, it follows from  (\ref{pde.1_1}) that
\be \label{cpde.1.2}g= d\rho^2+\sum_{i,j=1}^{n-1}\gamma_{ij}(r,\theta)d\theta_id\theta_j,\quad |x|\ge r_0,\ee
which implies $\rho(x)-c(r_0)$$(|x|\ge r_0)$ is the geodesic distance function of $(\R^n,g)$ from  $r_0(x/|x|)$ to $x$.

Let \be\label{cpde.1.1}\Upsilon(x)=(\gamma_{ij})_{(n-1)\times (n-1)}(x),\quad |x|\ge r_0.\ee

Let $S(r)$ be the sphere in $\R^n$ with
a radius $r$. We introduce a tensor field of rank 2 on $S(r)$ by
\be \label{cpde.1.2.1} \textbf{P}(X,Y)=\frac{1}{2\phi(r)}\sum_{i,j=1}^{n-1}\frac{\partial\gamma_{ij}}{\partial r}(r,\theta)X_i Y_j,\ee
where \be X=\sum_{i=1}^{n-1}X_i\frac{\partial}{\partial\theta_i}, Y=\sum_{j=1}^{n-1}Y_j\frac{\partial}{\partial\theta_j}\in S(r)_x,\quad |x|\ge r_0.\ee\\

The following Assumption ${\bf (B)}$ and Assumption ${\bf (C)}$ are the main assumptions.

{\bf Assumption (B)}\,\,\,  $g$ is a cone such that
\be  \phi(r)=m r^{m-1}, \ \ |x|\ge r_0,  \ee
where  $m>0$ is a constant.
\begin{rem} Let Assumption ${\bf (B)}$ hold true. It follows from  (\ref{cpde.2.2}) that \be \rho(x)=r^m,\quad  |x|\ge r_0. \ee
\end{rem}

If Assumption ${\bf (B)}$ holds true, we say that the cone g has {\bf polynomial of degree $m$ growth.}
$$ $$

{\bf Assumption (C)}\,\,\,  $g$ is a cone such that
\be \label{pde.3} \phi(r)=m r^{m-1}, \quad \textmd{for} \ \ |x|\ge r_0,  \ee
and
\be \label{pde.3-1} \textbf{P}(X,X)\ge \alpha(r) |X|^2_g, \quad \textmd{for} \ \ x\in \overline{\Om},X\in S(r)_x,\ee
where $m>0$ is a constant, $\alpha(r)$ is a positive function defined on $ [r_0,+\infty)$.
\begin {rem} If Assumption ${\bf (C)}$ holds true, Assumption ${\bf (B)}$ holds true.
\end {rem}
\begin{rem} Let Assumption ${\bf (C)}$ hold true. It follows from  (\ref{cpde.2.2}) that \be \rho(x)=r^n,\quad  |x|\ge r_0. \ee

\end{rem}

  The inequality (\ref{pde.3-1}) can be checked by the following proposition.
  \begin{pro} \label{1cpde.1}Let $g$ be a cone and $\alpha(r)$ be a smooth function defined on $[r_0,+\infty)$.
   Let
    \be\label{cpde.2} \Upsilon(r,\theta)= e^{\int_{r_0}^rh(y)dy}\Upsilon(r_0,\theta)+e^{\int_{r_0}^rh(y)dy}\int_{r_0}^r2\phi(y)e^{-\int_{r_0}^yh(z)dz} Q(y,\theta)dy,\quad |x|\ge r_0, \ee
    where $\Upsilon(r,\theta)$ is given by (\ref{cpde.1.1}) and

    \be h(r)= 2\alpha(r) \phi(r),\ee
      \be Q(x)=(q_{i,j})_{(n-1)\times (n-1)}(x),\quad q_{i,j}=\left\<\widehat{Q}(x)\frac{\partial}{\partial\theta_i},\frac{\partial}{\partial\theta_j}\right\>,\ee where $\widehat{Q}(x)$ is a smooth, symmetric, nonnegative definite matrix function defined on $|x|\ge r_0$.

 Then  \be  \label{cpde.1}{\rm{\textbf{P}}}(X,X)\ge \alpha(r) |X|^2_g, \quad \textmd{for }\ \ X\in S(r)_x,|x|\ge r_0.\ee
   \end{pro}

   {\bf Proof.}  \ \ \ It follows from (\ref{cpde.2}) that
    \be\frac{1}{2\phi(r)}\frac{\partial\Upsilon(x)}{\partial r}-\alpha(r) \Upsilon(x)=Q(x),\quad |x|\ge r_0.\ee
Since $\widehat{Q}(x)$ is a nonnegative definite matrix, we have
    \be\sum_{i,j=1}^{n-1}\left(\frac{1}{2\phi(r)}\frac{\partial\gamma_{ij}}{\partial r}-\alpha(r) \gamma_{ij}\right) X_i X_j= \left\<\widehat{Q}(x)X,X\right\> \ge 0, \quad \textmd{for} \ \  X\in S(r)_x,|x|\ge r_0.\ee
Then, the inequality (\ref{cpde.1}) holds true.$\Box$

Next, we give some examples of $A(x)$ which satisfy Assumption ${\bf (C)}$.
\begin{exl}
Let \be A(x)=\left(\frac{1}{\phi^2(r)}-e^{-\int^r_{r_0}hdy}\right)\frac{x\otimes x}{|x|^2}+e^{-\int^r_{r_0}hdy}I_n, \quad \textmd{for }\quad |x|\ge r_0,\ee
where $ \alpha(r)$ is a positive function defined on $ [r_0,+\infty)$ and
\be \phi(r)=mr^{m-1}, h(r)= 2(\alpha(r) -\frac{1}{r})\phi(r),\quad |x|\ge r_0.\ee

Then
\be A(x)x=\frac{1}{m^2r^{2(m-1)}}x  \quad \textmd{for} \ \ |x|\ge r_0,\ee
and
\be A(x)\frac{\partial}{\partial\theta_i}=e^{-\int^r_{r_0}hdy}\frac{\partial}{\partial\theta_i},\quad 1\leq i\leq n-1,|x|\ge r_0.\ee
Thus
\be \gamma_{ij}(r,\theta)= \left\<A^{-1}(r,\theta)\frac{\partial}{\partial\theta_i},\frac{\partial}{\partial\theta_j}\right\>=\left(\frac{r}{r_0}\right)^2e^{\int^r_{r_0}hdy}\gamma_{ij}(r_0,\theta),\quad 1\leq i,j\leq n-1,|x|\ge r_0,\ee
where
\be \gamma_{ij}(r_0,\theta)=\left\<\frac{\partial}{\partial\theta_i},\frac{\partial}{\partial\theta_j}\right\>(r_0,\theta).\ee

Therefore
    \be \Upsilon(r,\theta)= \left(\frac{r}{r_0}\right)^2e^{\int_{r_0}^rh(y)dy}\Upsilon(r_c,\theta).\ee
From Proposition \ref{1cpde.1}, we have
\be  {\rm{\textbf{P}}}(X,X)\ge \alpha(r) |X|^2_g \quad \textmd{for} \quad X\in S(r)_x,|x|\ge r_0.\ee

\end{exl}

\begin{exl}
Let \be A(x)=\frac{x\otimes x}{\phi^2(r)|x|^2}+e^{\int^r_{r_0}hdy}
Q\left(\frac{r_0x}{|x|}\right)\left(I_n-\frac{ x\otimes x}{|x|^2}\right)
, \quad \textmd{for}  \ \ |x|\ge r_0,\ee
where $ \alpha(r)$ is a positive function defined on $ [r_0,+\infty)$ and
\be \phi(r)=mr^{m-1}, h(r)= 2(\alpha(r) -\frac{1}{r})\phi(r),\quad |x|\ge r_0.\ee $Q(x)=(q_{i,j})_{n\times n}(x)$ is a symmetric,  definite matrix  for each $|x|= r_0$ and $q_{i,j}(x)(1\leq i,j \leq n)$ are smooth functions.

Then
\be A(x)x=\frac{1}{m^2r^{2(m-1)}}x  \quad for \ \ |x|\ge r_0,\ee
and
\be A(x)\frac{\partial}{\partial\theta_i}=e^{-\int^r_{r_0}hdy}Q\left(\frac{r_0x}{|x|}\right)\frac{\partial}{\partial\theta_i},\quad 1\leq i\leq n-1,|x|\ge r_0.\ee
Thus
\be \gamma_{ij}(r,\theta)= \left\<A^{-1}(r,\theta)\frac{\partial}{\partial\theta_i},\frac{\partial}{\partial\theta_j}\right\>=\left(\frac{r}{r_0}\right)^2e^{\int^r_{r_0}hdy}\gamma_{ij}(r_0,\theta),\quad 1\leq i,j\leq n-1,|x|\ge r_0,\ee
where
\be \gamma_{ij}(r_0,\theta)=\left\<Q^{-1}\frac{\partial}{\partial\theta_i},\frac{\partial}{\partial\theta_j}\right\>(r_0,\theta).\ee

Therefore
    \be \Upsilon(r,\theta)= \left(\frac{r}{r_0}\right)^2e^{\int_{r_0}^rh(y)dy}\Upsilon(r_0,\theta).\ee
From Proposition \ref{1cpde.1}, we have
\be  {\rm{\textbf{P}}}(X,X)\ge \alpha(r) |X|^2_g, \quad for \ \  X\in S(r)_x,|x|\ge r_0.\ee

\end{exl}

\vskip .2cm
\subsection{Main Theorems}
  \quad \ \ In order to facilitate the discussion, we define
 \be \nu(x) \ \ is\ \  the\ \   unit\ \  normal \ \ vector\ \ outside\ \  \Om  \ \ in \ \ \R^n\ \  for\ \  x\in  \Ga, \ee
    \be \label{pde.3.1}\nu_{\A}(x)=A(x)\nu,\quad  x\in \partial \Om,  \ \ and  \  \ \nu_{\A}(x)=A(x)\frac{\partial}{\partial r}, \quad |x|> r_1.\ee
Note that \be \label{1pde.1} |\nu_{\A}|_g=\sqrt{\<A(x) \frac{\partial }{\partial r},A(x) \frac{\partial }{\partial r}\>_g}=\sqrt{\< \frac{\partial }{\partial r},A(x) \frac{\partial }{\partial r}\>}\le F(|x|),\quad |x|>r_1,\ee
where  $F(|x|)$ is given by (\ref{wg.4_1}).
$$ $$

Our two primary decay results are now listed as follows:
\begin{thm}\label{wg.7} Suppose that the following three conditions hold:
\be \label{wg.7_1}Assumption\ {\bf(C)}\ holds\ true\ with \ m=n\ and   \ \alpha(r)= \delta\frac{e^{\delta r^n }+e^{-\delta r^n }}{e^{\delta r^n }-e^{-\delta  r^n }}, \quad |x|\ge r_0, \ee
\be \label{wg.7_2}\frac{\partial r}{\partial \nu}\le 0,\quad x\in \Ga,\ee
\be \label{wg.7_4}u_0(x)=u_1(x)=0,\quad |x|\ge R_0, \ee
where $\delta>0,R_0 >r_1$ are  constants.

Then
there
exist positive constants $C_1(a,R_0),\ C_2
>0$ such that
\begin{equation}
\label{wg.7.1}
E(a,t)\leq C_1(a,R_0) e^{-C_2 t}E(0),\quad \forall t> 0.\end{equation}
\end{thm}
\begin{rem}
The estimate (\ref{wg.7.1}) is independent of the parity of the dimension $n$, which is the main difference with the constant coefficient wave equation.\end{rem}
\begin{rem}
For $n\ge 2$, the condition (\ref{wg.7_1}) does not hold at the original point $O$. Thus, the above results may not hold for the wave equation on $\R^n$. \end{rem}

\begin{thm}\label{wg.7.2_1}  Suppose that the following three conditions hold:
\be \label{wg.7.2_3}Assumption\ {\bf (C)} \ \ holds\ \ true\ \ with  \ \  \alpha(r)=  \frac{m_1}{ r^m},\quad |x|\ge r_0, \ee
\be \frac{\partial r}{\partial \nu}\le 0,\quad x\in \Ga,\ee
\be \label{wg.39}u_0(x)=u_1(x)=0,\quad |x|\ge R_0, \ee
where $m_1>\frac{1}{2}, R_0>r_1$ are constants.

Then there exists a positive constant $C(a,R_0)>0$ such that
\begin{equation}
\label{wg.7.2_2}
E(a,t)\leq\frac{ C(a,R_0)}{t} E(0),~~\forall~ t> 0.
\end{equation}
\end{thm}

 \vskip .5cm
\def\theequation{3.\arabic{equation}}
\setcounter{equation}{0}
\section{Multiplier Identities and Key Lammas}
\vskip .2cm

\quad \ \ We need to establish several multiplier identities, which are useful for our problem.

\begin{lem}\label{wg.14}
Suppose that $u(x,t)$ solves the system (\ref{wg.1})
and  $\H(x,t)=\Sigma_{i=1}^n h_i(x,t)\frac{\partial}{\partial x_i}$ is a  time-varying vector
field defined on $\overline \Om$, where $h_i\in C^1(\overline{\Om}\times [0,+\infty))$. Then
\begin{eqnarray}
 \label{wg.14.1}
 &&\int_0^T\int_{\partial\Om(a)}\frac{\pa u}{\pa\nu_{\A}}\H(u) d\Ga dt+\frac12\int_0^T\int_{\partial\Om(a)}
\left(u_t^2-\left|\nabla_g u\right|_g^2\right)\H\cdot\nu d\Ga dt\nonumber\\
=&&(u_t,\H(u))\Big |^T_0+\int_0^T\int_{\Om(a)}D\H(\nabla_g
u,\nabla_g u) dx dt-\int_0^T\int_{\Om(a)}u_t \H_t(u)dxdt\nonumber\\
\quad &&+\frac12\int_0^T\int_{\Om(a)}\left(u_t^2-\left |\nabla_g
u\right|_g^2\right)\div\H dx dt,
\end{eqnarray}
where
\be (u_t,\H(u))\Big |^T_0 =\int_{\Omega(a)}u_t\H(u) dx\Big |^T_0.\ee

Moreover, assume that $P\in C^1(\overline{\Om}\times [0,+\infty))$. Then
\begin{eqnarray}
\label{wg.14.2}
\int_0^T\int_{\Om(a)}\left(u_t^2-\left |\nabla_g
u\right|_g^2\right)P dx dt =&& \frac12\int_0^T
\int_{\Om(a)}\nabla_gP(u^2) dx dt-\int_0^T\int_{\partial\Om(a)}Pu\frac{\pa
u}{\pa\nu_{\A}}d\Ga dt\nonumber\\
&&+(u_t,u P)\Big |^T_0-\int_0^T
\int_{\Om(a)}u_tu P_tdxdt,\end{eqnarray}
and
\begin{eqnarray}
\label{wg.14.3}
 \int_{\Om(a)}\left(u_t^2+\left |\nabla_g
u\right|_g^2\right)P dx \Big|^T_0=&&\int_0^T
\int_{\Om(a)}P_t\left(u_t^2+\left |\nabla_g
u\right|_g^2\right) dxdt\nonumber\\&&-2\int_0^T
\int_{\Om(a)}u_t\nabla_gP(u) dx dt\nonumber\\
&&+2\int_0^T\int_{\partial\Om(a)}Pu_t\frac{\pa
u}{\pa\nu_{\A}}d\Ga dt.\end{eqnarray}

\end{lem}

{\bf Proof}.
 Firstly, we multiply the wave equation in (\ref{wg.1}) by $\H(u)$ and integrate over $\Omega(a)\times
(0,T) $, noting that
\begin{eqnarray} \<\nabla_g u,\nabla_g (\H(u))\>_g =&&\nabla_g u\<\nabla_g u,\H\>_g
= D^2u(\H,\nabla_g u)+D\H(\nabla_g u,\nabla_g u)\nonumber\\
=&& D^2u(\nabla_g u,\H)+D\H(\nabla_g u,\nabla_g u)\nonumber\\
=&&\frac{1}{2}\H(|\nabla_g u|_g^2)+D\H(\nabla_g u,\nabla_g u)\nonumber\\
=&&D\H(\nabla_g u,\nabla_g u)+\frac{1}{2}\div(|\nabla_g u|_g^2\H)-\frac{1}{2}|\nabla_g u|_g^2\div \H.\end{eqnarray}
the equality  (\ref{wg.14.1}) follows from Green's formula.

Secondly, we multiply the wave equation in (\ref{wg.1}) by $Pu$ and integrate over $\Omega(a)\times
(0,T)$. The equality  (\ref{wg.14.2}) follows from Green's formula.
Finally, the equality  (\ref{wg.14.3}) follows from Green's formula.$\Box$

The following three lemmas will be utilized frequently in our subsequent proof.
 \begin{lem}\label{pde.4} Let $g$ be a cone. Then
 \be \label{pde.5} D^2\rho(X,X)= {\rm{\textbf{P}}}(X,X),\quad \textmd{for}\ \  X\in S(r)_x, |x|\ge r_0, \ee
 where $\rho(x)$ is given by (\ref{cpde.2.1}), ${\rm{\textbf{P}}}(.,.)$ is given by  (\ref{cpde.1.2.1}) and $D^2\rho$ is the Hessian of $\rho$ in the metric $g$.
  \end{lem}

  {\bf Proof}. Let $\Upsilon(x)$ be given by (\ref{cpde.1.1}). Denote
  \be \aleph=\left\{x\ \ \Big|\quad |x|\ge r_0, \det(\Upsilon(x))\neq 0\right\}.\ee
 Let $x\in \aleph$ and denote $\theta_n=r$, we have
   \be \gamma_{ni}(x)=\gamma_{in}(x)=\left\<A^{-1}(x)\frac{\partial}{\partial \theta_n},\frac{\partial}{\partial \theta_i}\right\>=0\ \ 1\le i\le n-1,\quad \gamma_{nn}(x)=1.\ee
 Let  $(\gamma^{ij})_{n\times n}(x)=(\gamma_{ij})^{-1}_{n\times n}(x)$.  Then \be \gamma^{ni}(x)=\gamma^{in}(x)=0\ \ 1\le i\le n-1,\quad \gamma^{nn}(x)=1.\ee

 We compute Christofell symbols as
  \be \Gamma^k_{in}=\frac{1}{2}\sum_{l=1}^n\gamma^{kl}\left(\frac{\partial(\gamma_{il})}{\partial r}+\frac{\partial(\gamma_{nl})}{\partial\theta_i}-\frac{\partial(\gamma_{in})}{\partial\theta_l}\right)
  =\frac{1}{2}\sum_{l=1}^{n-1}\gamma^{kl}\frac{\partial(\gamma_{il})}{\partial r},\ee
  for $ 1\le i,k\le n-1$,  which give
  \be D_{\frac{\partial}{\partial \theta_i}}\frac{\partial}{\partial r}=\frac{1}{2}\sum_{k=1}^{n-1}\left(\sum_{l=1}^{n-1}\gamma^{kl}\frac{\partial(\gamma_{il})}{\partial r}\right)\frac{\partial}{\partial \theta_k}.\ee
Then for $X=\sum_{i=1}^{n-1}X_i\frac{\partial}{\partial\theta_i} \in S(r)_x$, we deduce that
\begin{eqnarray}  D^2\rho(X,X)&&= \sum_{i,j=1}^{n-1}\left\< D_{\frac{\partial}{\partial \theta_i}}\frac{\partial}{\partial \rho},\frac{\partial}{\partial \theta_j}\right\>_gX_iX_j
  \nonumber\\
&& =\frac{1}{2}\sum_{i,j,k,l=1}^{n-1}\gamma^{kl}\frac{\partial(\gamma_{il})}{\partial\rho}\gamma_{kj}X_iX_j \nonumber\\
&& =\frac{1}{2}\sum_{i,j=1}^{n-1}\frac{\partial(\gamma_{ij})}{\partial\rho}X_iX_j  =\frac{1}{2\phi(r)}\sum_{i,j=1}^{n-1}\frac{\partial(\gamma_{ij})}{\partial r}X_iX_j={\rm{\textbf{P}}}(X,X).\end{eqnarray}

Note that
\be \aleph \ \  is\ \  dense\ \  in \ \ |x|\ge r_0. \ee
The equality (\ref{pde.5}) holds true. $\Box$

\begin{lem}\label{wg.8} Let Assumption ${\bf (A)}$ hold true. Let $\psi(y)$ be a nonnegative function defined on $[r_0,+\infty)$ satisfying
\be \label{wg.8.1}\int_{r_0}^{+\infty}\psi(y)dy<+\infty. \ee
Then
      \be \label{wg.8.2}\lim_{\overline{y\rightarrow +\infty}}F(y)\psi(y)=0, \ee  where  $F(y)$ is given by (\ref{wg.4_1}). \end{lem}

{\bf Proof.\ }If (\ref{wg.8.2}) does not hold true, then there exist constants $\sigma>0, M>r_0$ satisfying
\be F(y)\psi(y)\ge \sigma,  \quad  \forall y\ge M. \ee
Then
   \be \psi(y)\ge \frac{\sigma}{F(y)},  \quad \forall  y\ge M. \ee
With (\ref{wg.4}),  we have
   \be \int_{r_0}^{+\infty}\psi(y)dy\ge \int_{M}^{+\infty}\psi(y)dy\ge\int_{M}^{+\infty}\frac{\sigma}{F(y)} dy =+\infty, \ee
which contradicts (\ref{wg.8.1}).\hfill $\Box$

\begin{lem}\label{wg.9}Let Assumption ${\bf (A)}$ hold true. Let $u,v$ satisfy
   \be \label{wg.9.2}\int_\Om (|\nabla_g u|^2_g+v^2)dx<+\infty. \ee
Then
      \be \label{wg.9.1}\lim_{\overline{y\rightarrow +\infty}}\int_{x\in \Om,|x|=y} \left|v\frac{\pa u}{\pa\nu_{\A}}\right|=0. \ee \end{lem}
{\bf Proof.\ }  With (\ref{pde.3.1}) and (\ref{1pde.1}), we deduce that
 \be \int_\Om (|\nabla_g u|^2_g+v^2)dx=\int_{r_0}^{+\infty}dy\int_{x\in \Om,|x|=y}  (|\nabla_g u|^2_g+v^2)dx,\ee
and for $y>r_1$
   \begin{eqnarray} \left|v\frac{\pa u}{\pa\nu_{\A}}\right|&&= |v\<\nabla_g u, \nu_{\A}\>_g |\nonumber\\
&&\le |v||\nabla_g u|_g |\nu_{\A}|_g\nonumber\\
&&\le |v||\nabla_g u|_g F(y)
\nonumber\\
&&\le \frac{ 1}{2}(|\nabla_g u|^2_g+v^2)F(y),\end{eqnarray}
the estimate (\ref{wg.9.1}) follows from Lemma \ref{wg.8}.\hfill $\Box$
 \vskip .5cm
\def\theequation{4.\arabic{equation}}
\setcounter{equation}{0}
\section{Proofs for   well-posedness and   propagation property }
\vskip .2cm

\quad \ \ {\bf Proof of Theorem \ref{wg.5} }\,\,\,
Let
$$\A=\left(\matrix{0&1\cr
\div A(x)\nabla &0\cr}\right).$$
Then the system (\ref{wg.1}) can be rewritten as
\begin{equation}
 \cases{ U_t-\A U=0\qquad (x,t)\in \Omega\times
(0,+\infty),\cr
 U\large |_{\Ga}=0\qquad t\in(0,+\infty),\cr
U(0)=U_0\qquad x\in \Omega,}
\end{equation}
where
$$U=\left(\matrix{u\cr
v\cr}\right), U_0=\left(\matrix{u_0\cr
u_1\cr}\right).$$

Define $V$ and $\<,\>_V$ by
\be  V=\widehat{H}^1_{0}(\Om)\times L^2(\Om), \ee
\be \<U_1,U_2\>_{V}= \int_{\Omega}\<\nabla_g u_1,\nabla_g u_2\>_g dx+\int_\Om v_1v_2dx,\ee
where $$U_1=\left(\matrix{u_1\cr
v_1\cr}\right), U_2=\left(\matrix{u_2\cr
v_2\cr}\right).$$
Note that
 \be D(\A)=\widehat{H}^2_{0}(\Om)\times  (\widehat{H}^1_{0}(\Om)\cap  L^2(\Om)).\ee

 The next is to prove the operator $i\A$ is self-joint.

Let $U_1,U_2\in D(\A)$. With (\ref{wg.9.1}), we deduce that
 \begin{eqnarray}\<\A U_1,U_2\>_{V}&&=\int_{\Omega}\<\nabla_g v_1,\nabla_g u_2\>_g dx+\int_\Om v_2\div A(x)\nabla u_1dx\nonumber\\
&&=\int_{\Omega}\<\nabla_g v_1,\nabla_g u_2\>_g dx-\int_{\Omega}\<\nabla_g u_1,\nabla_g v_2\>_g dx,
\end{eqnarray}
 and
 \begin{eqnarray}\< U_1,\A U_2\>_{V}&&=\int_{\Omega}\<\nabla_g u_1,\nabla_g v_2\>_g dx+\int_\Om v_1\div A(x)\nabla u_2dx\nonumber\\
&&=\int_{\Omega}\<\nabla_g u_1,\nabla_g v_2\>_g dx-\int_{\Omega}\<\nabla_g v_1,\nabla_g u_2\>_g dx.
\end{eqnarray}
Then $ \<\A U_1,U_2\>_{V}=- \< U_1,\A U_2\>_{V}$, the operator $i\A$ is symmetric.

If \be(i\A\pm iI)U=0,\ee then
\be\cases{u\pm v=0,\cr
v\pm \div A(x)\nabla u=0.}\ee
We obtain $u=v=0$, that is $U=0$. Thus, the operator $i\A$ is self-joint.

 By Stone theorem (Theorem 10.8 in \cite{w17}), we have $U(t)= e^{t\A}U(0)$.\hfill $\Box$\\
 $$ $$

{\bf Proof of Theorem \ref{wg.c5} }\,\,\,
Let
\be  \widetilde{E}(t)=\frac12\int_{\Omega\setminus\Omega(z +R_0)}\left(u_t^2+|\nabla_gu|^2_g\right)dx.\ee
With  (\ref{pde.3.1}) and (\ref{1pde.1}), we deduce that
\begin{eqnarray}
\widetilde{E}'(t) && =\frac12\frac{d}{dt}\int_{z +R_0}^{+\infty}dy\int_{x\in \Om,|x|=y} (u_t^2+|\nabla_gu|^2_g)dx\nonumber\\
&&=-\frac12\int_{|x|=z +R_0}\frac{dz}{dt}\left(u_t^2+|\nabla_gu|^2_g\right)d\Ga +\int_{\Omega\setminus\Omega(z +R_0)}\div u_t \nabla_g udx
 \nonumber\\
&&=-\frac12\int_{|x|=z +R_0}F(z +R_0)\left(u_t^2+|\nabla_gu|^2_g\right)d\Ga -\int_{|x|=z +R_0} u_t\frac{\pa u}{\pa\nu_{\A}} d\Ga
\nonumber\\
&&\le-\frac12\int_{|x|=z +R_0}F(z +R_0)\left(u_t^2+|\nabla_gu|^2_g\right)d\Ga  +\left|\int_{|x|=z +R_0} u_t|\nabla_gu|_g\cdot |\nu_{\A}|_g d\Ga\right|\nonumber\\
&&\le-\frac12\int_{|x|=z +R_0}F(z +R_0)\left(u_t^2+|\nabla_gu|^2_g\right)d\Ga  +\frac12\int_{|x|=z +R_0} F(z +R_0)\left(u_t^2+|\nabla_gu|^2_g\right)d\Ga\nonumber\\
&&=0.
\end{eqnarray}
Noting that $\widetilde{E}(0)=0$, we have $u(x,t)=0,\quad |x|\ge z +R_0$.\hfill $\Box$

 \vskip .5cm
\def\theequation{5.\arabic{equation}}
\setcounter{equation}{0}
\section{Uniform decay of local energy for radial solutions }
\quad \ In this section, we shall study the differences of the decay rate of the local energy between the variable coefficients wave equation and  the   constant coefficient wave equation for radial solutions.

    Let $\Ga=\{x|\ |x|=r_0\}$ and let $u_0(x), u_1(x)$ be of compact support.

It is well-known  that the local energy for the constant coefficient wave equation
\begin{equation}
\label{wg.c11} \cases{u_{tt}- \Delta  u=0\qquad (x,t)\in \Om\times
(0,+\infty),\cr
 u|_{\Ga}=0\qquad t\in(0,+\infty),
\cr u(x,0)=u_0(x),\quad u_t(x,0)=u_1(x)\qquad x\in \Omega,}
\end{equation}
has a uniform decay rate as follows. For  even dimensional space, the uniform decay rate of the local energy is polynomial and  for  odd dimensional space, the uniform decay rate of the local energy  is exponential. See \cite{w7}, \cite{w4}.

Let Assumption ({\bf B}) hold true.  From (\ref{cpde.2.2}), we have \be \label{cpde.2.3} \rho(x)=r^m,\quad x\in \Om.\ee

Let $u_0(x)=u_0(\rho), u_1(x)=u_1(\rho)$
 and let $u(\rho,t)$ solve the following system
\begin{equation}
\label{wg.11} \cases{u_{tt}-u_{\rho\rho}-(\frac{n}{m}-1)\frac{1}{\rho}u_{\rho}=0\qquad (\rho,t)\in (r_0^m,+\infty)\times
(0,+\infty),\cr
 u(r_0^m)=0,
\cr u(0)=u_0(\rho),\quad u_t(0)=u_1(\rho)\qquad \rho\ge r_0^m.}
\end{equation}

Note that \be \label{wg.11.1}
\div \frac{\partial}{\partial \rho}=\div \frac{dr}{d\rho}\frac{\partial}{\partial r}=\frac{1}{m}\div r^{1-m}\frac{\partial}{\partial r}=\frac{n-m}{m}\frac{1}{r^m}=(\frac{n}{m}-1)\frac{1}{\rho},\quad \textmd{for} \ \ x \in \Om.\ee
Then
\begin{eqnarray} &&\label{wg.11_1}\div A(x) \nabla u =\div(u_{\rho}\frac{\partial}{\partial \rho}+\nabla_{\Ga_g}u)=u_{\rho\rho}+(\frac{n}{m}-1)\frac{1}{\rho}u_{\rho}, \quad \textmd{for} \ \ x \in \Om,\end{eqnarray}
where
\be \nabla_{\Ga_g}u=\nabla_gu-u_\rho\frac{\partial}{\partial\rho}.\ee
Thus, $u(x,t)=u(\rho(x),t)$ solves the  system (\ref{wg.1}). Then the energy and the local energy for the system (\ref{wg.1})  can be rewritten as

 \be E(t,a)=\frac12\int_{\Omega(a)}\left(u_t^2+u_\rho^2\right)dx,\ee
\begin{equation}
 E(t)=\frac12\int_{\Omega}\left(u_t^2+u_\rho^2\right)dx.
\end{equation}

Using the conclusion of  scattering theory of the constant coefficient wave equation, for any positive integer $k$ , we have

\begin{itemize}
\item
If $\frac{n}{m}=2k-1$ in (\ref{wg.11}), which implies $m=\frac{n}{2k-1}$, the   decay rate of the local energy for the system (\ref{wg.1}) is exponential, whether the dimension $n$ is even or odd.
\item If $\frac{n}{m}=2k$ in (\ref{wg.11}), which implies $m=\frac{n}{2k}$,  the   decay rate of the local energy for the system (\ref{wg.1}) is polynomial, whether the dimension $n$ is even or odd.
    \end{itemize}

\vskip .5cm
\def\theequation{6.\arabic{equation}}
\setcounter{equation}{0}
\section{Uniform decay of local energy for general solutions}
\vskip .2cm
\subsection{Exponential  decay  of local energy }

$$ $$

 {\bf Proof of Theorem \ref{wg.7} }\,\,\,

Let $\rho(x)$ be given by (\ref{cpde.2.2}), then
\be  \rho(x)=r^n, x\in\overline{\Om}. \ee

 With (\ref{wg.7_1}) and (\ref{pde.5}), we have
\be D^2\rho(X,X) \ge \delta \frac{e^{\delta \rho }+e^{-\delta \rho }}{e^{\delta \rho }-e^{-\delta \rho }}|X|_g^2,\quad \textmd{for}\ \  X\in S(r)_x, x\in\Om. \ee
Let $\H=e^{\delta t}(e^{\delta\rho}-e^{-\delta\rho})\frac{\partial}{\partial \rho}$ in (\ref{wg.14.1}). For $x\in \Om$, we deduce that
\begin{eqnarray} \label{wg1.16_1}D\H(\nabla_gu,\nabla_g u)&&=D\H(u_\rho\frac{\partial}{\partial \rho} ,u_\rho\frac{\partial}{\partial \rho})+D\H(\nabla_{\Ga_g}u,\nabla_{\Ga_g}u)\nonumber\\
&&=\delta e^{\delta t}(e^{\delta\rho}+e^{-\delta\rho})u^2_\rho+e^{\delta t}(e^{\delta\rho}-e^{-\delta\rho})D^2\rho(\nabla_{\Ga_g}u,\nabla_{\Ga_g}u)
\nonumber\\
&&\ge\delta e^{\delta t}(e^{\delta\rho}+e^{-\delta\rho})(u^2_\rho+ |\nabla_{\Ga_g}u|^2_g)=\delta e^{\delta t}(e^{\delta\rho}+e^{-\delta\rho})|\nabla_{g}u|^2_g,\end{eqnarray}
where
$\nabla_{\Ga_g}u=\nabla_gu-u_\rho\frac{\partial}{\partial\rho}$.
Note that \be
\div \frac{\partial}{\partial \rho}=\div \frac{dr}{d\rho}\frac{\partial}{\partial r}=\frac{1}{n}\div r^{1-n}\frac{\partial}{\partial r}=0,\quad x \in \Om.\ee
Then
\begin{eqnarray}\label{wg1.16_2}\div \H&&=\delta e^{\delta t}(e^{\delta\rho}+e^{-\delta\rho})+e^{\delta t}(e^{\delta\rho}-e^{-\delta\rho})\div\frac{\partial}{\partial \rho} \nonumber\\
&&=\delta e^{\delta t}(e^{\delta\rho}+e^{-\delta\rho}),\quad x \in \Om.\end{eqnarray}

 Substituting (\ref{wg1.16_1}) and (\ref{wg1.16_2}) into (\ref{wg.14.1}), we have
\begin{eqnarray}
 \label{wg1.16}
 &&\int_0^T\int_{\partial\Om(a)}\frac{\pa u}{\pa\nu_{\A}}\H(u) d\Ga dt+\frac12\int_0^T\int_{\partial\Om(a)}
\left(u_t^2-\left|\nabla_g u\right|_g^2\right)\H\cdot\nu d\Ga dt\nonumber\\
\ \ &&\quad\geq(u_t,\H(u))\Big |^T_0+\int_0^T\int_{\Om(a)}\delta e^{\delta t}(e^{\delta\rho}+e^{-\delta\rho})|\nabla_{g}u|^2_g dx dt\nonumber\\
\ \ &&\qquad-\int_0^T\int_{\Om(a)}\delta e^{\delta t}(e^{\delta\rho}-e^{-\delta\rho})u_t u_\rho dx dt +\int_0^T\int_{\Om(a)}\frac{\delta e^{\delta t}(e^{\delta\rho}+e^{-\delta\rho})}{2}\left(u_t^2-\left |\nabla_g
u\right|_g^2\right) dx dt\nonumber\\
\ \ &&\quad=(u_t,\H(u))\Big |^T_0-\int_0^T\int_{\Om(a)}\delta e^{\delta t}(e^{\delta\rho}-e^{-\delta\rho})u_t u_\rho dx dt\nonumber\\
\ \ &&\qquad+\int_0^T\frac{\delta e^{\delta t}(e^{\delta\rho}+e^{-\delta\rho})}{2}\int_{\Om(a)}\left(u_t^2+\left |\nabla_g
u\right|_g^2\right) dx dt.
\end{eqnarray}
Let $P=\frac{e^{\delta t}(e^{\delta\rho}+e^{-\delta\rho})}{2}$, substituting (\ref{wg.14.3}) into (\ref{wg1.16}),  letting $a\rightarrow+\infty$,  we obtain
 \begin{eqnarray}
 \label{wg1.16.1}
 &&\Pi_{\Ga}\geq(u_t,\H(u))\Big |^T_0+\frac{1}{2}\int_{\Om(a)} e^{\delta t}(e^{\delta\rho}+e^{-\delta\rho})\left(u_t^2+\left |\nabla_g
u\right|_g^2\right) dx \Big|^T_0,
\end{eqnarray}
where
\begin{eqnarray}
 \label{wg1.16.1.111}
\Pi_{\Ga}=&&\int_0^T\int_{\Ga}\frac{\pa u}{\pa\nu_{\A}}\H(u) d\Ga dt+\frac12\int_0^T\int_{\Ga}
\left(u_t^2-\left|\nabla_g u\right|_g^2\right)\H\cdot\nu d\Ga dt\nonumber\\
&&+2\int_0^T\int_{\Ga}Pu_t\frac{\pa
u}{\pa\nu_{\A}}d\Ga dt.
\end{eqnarray}

 Since $u\Large|_{\Ga}=0,$ we
obtain $\nabla_{\Ga_g} u\large|_{\Ga}=0$, that is,
 \be\nabla_g u=\frac{\pa
 u}{\pa\nu_{\A}}\frac{\nu_{\A}}{|\nu_{\A}|_g^2}\quad\mbox{for}\quad
 x\in\Ga.\label{wg1.16.1.2}\ee
Similarly, we have \be \H(u)=\<\H,\nabla_g u\>_g=\frac{\pa
 u}{\pa\nu_{\A}}\frac{\H\cdot\nu}{|\nu_{\A}|_g^2}\quad\mbox{for}\quad
 x\in\Ga.\label{wg1.16.1.3}\ee
With (\ref{wg.7_1}) and  (\ref{wg.7_2}), we deduce that
\be\label{wg1.16.1_2} \left\<\frac{\partial}{\partial \rho},\nu\right\>=\left\<\frac{\partial}{\partial r}\frac{dr}{ d \rho},\nu\right\>=\frac{r^{1-n}}{n} \left\<\frac{\partial}{\partial r},\nu\right\>=\frac{r^{1-n}}{n} \frac{\partial r}{\partial \nu}\le 0 \quad \mbox{for} \ \ x\in \Ga.\ee
Using the formula (\ref{wg1.16.1.3}) in the formula
(\ref{wg1.16.1.111}) on the portion $\Ga$, with (\ref{wg1.16.1_2}), we obtain
\begin{equation}
\label{wg1.16.1.4} \Pi_{\Ga}=\frac12\int_0^T
\int_{\Ga}\left(\frac{\pa
u}{\pa\nu_{\A}}\right)^2\frac{\H\cdot\nu}{|\nu_{\A}|_g^2} d\Ga dt\le 0.
\end{equation}

Substituting (\ref{wg1.16.1.4}) into (\ref{wg1.16.1}), we have
 \begin{eqnarray}
 \label{wg1.16.2}
 &&(u_t,\H(u))\Big |^T_0+\frac{1}{2}\int_{\Om} e^{\delta t}(e^{\delta\rho}+e^{-\delta\rho})\left(u_t^2+\left |\nabla_g
u\right|_g^2\right) dx \Big|^T_0\le0.
\end{eqnarray}

Accordingly,
 \begin{eqnarray}
 \label{wg1.16.3}
 &&\int_{\Om} e^{\delta T}(e^{\delta\rho}+e^{-\delta\rho})\left(u_t^2+\left |\nabla_g
u\right|_g^2\right) dx \nonumber\\
&&\qquad \leq\int_{\Om} e^{\delta T}(e^{\delta\rho}-e^{-\delta\rho})\left(u_t^2+u_\rho^2\right) dx+C(R_0)E(0) .
\end{eqnarray}
Then
 \begin{eqnarray}
 \label{wg1.16.4}
 &&2\int_{\Om} e^{\delta (T-\rho)}\left(u_t^2+u_\rho^2\right) dx+ \int_{\Om} e^{\delta T}(e^{\delta\rho}+e^{-\delta\rho})\left |\nabla_{\Ga_g}
u\right|_g^2 dx\leq C(R_0)E(0) .
\end{eqnarray}
Therefore
\begin{eqnarray}
 &&\int_{\Om} e^{\delta (T-\rho)}\left(u_t^2+\left |\nabla_g
u\right|_g^2\right) dx\leq C(R_0)E(0) .
\end{eqnarray}
The estimate (\ref{wg.7.1}) holds.  \hfill $\Box$

\vskip .2cm
\subsection{Polynomial decay  of local energy }
\begin{lem}\label{wg.17} Let Assumption ${\bf (B)}$ hold true and $u_0,u_1$ be of compact support. Let u solve the system (\ref{wg.1}).
Then
  \be \label{wg.17.1} \int_\Om u u_\rho dx=\int_\Om \frac{1}{2\rho}\left(1- \frac{n}{m}\right)u^2dx ,\quad \forall t\ge 0,\ee
  where $\rho(x)$ is given by (\ref{cpde.2.2}).
\end{lem}

{\bf Proof}. Note that
\be\label{wg.17_1}\div u^2\frac{\partial}{\partial \rho} =u^2 \div\frac{\partial}{\partial \rho} +2uu_\rho=\frac{n-m}{m\rho}u^2 +2uu_\rho. \ee
Integrate (\ref{wg.17_1}) over $\Om$,
the equality (\ref{wg.17.1}) holds.$\Box$

\begin{lem}\label{2wg.3}  Let Assumption ${\bf (B)}$ hold true and  $u_0,u_1$ be of compact support.  Let u solve the system (\ref{wg.1}).
 Then
\be \label{1wg.20}\int_{\Om}\rho\left(u_t^2+\left |\nabla_g
u\right|_g^2\right) dx \leq \int_0^T\int_{\Om} (u_t^2+u_\rho^2) dx dt+ \int_{\Om} \rho\left(u_1^2+\left |\nabla_g
u_0\right|_g^2\right) dx,\quad \forall t\ge 0,  \ee
 where $\rho(x)$ is given by (\ref{cpde.2.2}).
\end{lem}

{\bf Proof}.
Let $Q=(\rho-t)$, letting $a\rightarrow+\infty$, it follows from (\ref{wg.14.3}) that
\begin{eqnarray}
 \int_{\Om}(\rho-t)\left(u_t^2+\left |\nabla_g
u\right|_g^2\right) dx \Big|^T_0
=&& -\int_0^T
\int_{\Om}\left(u_t^2+\left |\nabla_g
u\right|_g^2\right)  dxdt-2\int_0^T
\int_{\Om}u_tu_\rho dx dt\nonumber\\
\leq&&-\int_0^T
\int_{\Om}\left |\nabla_{\Ga_g}
u\right|_g^2 dxdt.\end{eqnarray}
Then
\begin{eqnarray}
 \int_{\Om}(\rho-T)\left(u_t^2+\left |\nabla_g
u\right|_g^2\right) dx &&\leq-\int_0^T
\int_{\Om}\left |\nabla_{\Ga_g}
u\right|_g^2 dxdt\nonumber\\
&&\quad + \int_{\Om} \rho\left(u_1^2+\left |\nabla_g
u_0\right|_g^2\right) dx.\end{eqnarray}
Simple calculation shows that
\begin{eqnarray}
 \int_{\Om}\rho\left(u_t^2+\left |\nabla_g
u\right|_g^2\right)  dx&&\leq TE(0)  -\int_0^T
\int_{\Om}\left |\nabla_{\Ga_g}
u\right|_g^2 dx dt + \int_{\Om} \rho\left(u_1^2+\left |\nabla_g
u_0\right|_g^2\right) dx
\nonumber\\
&&= \int_0^T\int_{\Om} (u_t^2+u_\rho^2) dxdt+ \int_{\Om} \rho\left(u_1^2+\left |\nabla_g
u_0\right|_g^2\right) dx.\end{eqnarray}
The estimate (\ref{1wg.20}) holds.
$\Box$

\begin{lem}\label{2wg.5}  Let Assumption ${\bf (B)}$ hold true  and   $u_0,u_1$ be of compact support.  Let u solve the system (\ref{wg.1}).
 Then
\be \label{1wg.20.3} \int_{\Om}e^{\rho-t}\left(u_t^2+\left |\nabla_g
u\right|_g^2\right) dx \leq \int_{\Om}e^{\rho}\left(u_1^2+\left |\nabla_g
u_0\right|_g^2\right)  dx,\quad \forall t\ge 0,  \ee
 where $\rho(x)$ is given by (\ref{cpde.2.2}).
\end{lem}

{\bf Proof}.
Let $Q=e^{\rho-t}$, letting $a\rightarrow+\infty$, it follows from (\ref{wg.14.3}) that
\begin{eqnarray}
 \int_{\Om}e^{\rho-t}\left(u_t^2+\left |\nabla_g
u\right|_g^2\right) dx \Big|^T_0=&&-\int_0^T
\int_{\Om}e^{\rho-t}\left(u_t^2+\left |\nabla_g
u\right|_g^2\right)  dxdt\nonumber\\&&-2\int_0^T
\int_{\Om}e^{\rho-t}u_tu_\rho dx dt\nonumber\\
\leq&&-\int_0^T
\int_{\Om}e^{\rho-t}\left(u_t^2+\left |\nabla_g
u\right|_g^2\right)  dxdt\nonumber\\&&+\int_0^T
\int_{\Om}e^{\rho-t}(u^2_t+u^2_\rho )dx dt\nonumber\\
\leq&& 0.\end{eqnarray}
Then
\be   \int_{\Om}e^{\rho-T}\left(u_t^2+\left |\nabla_g
u\right|_g^2\right) dx  \leq \int_{\Om}e^{\rho}\left(u_1^2+\left |\nabla_g
u_0\right|_g^2\right) dx.\ee
The estimate (\ref{1wg.20.3}) holds.
$\Box$

\begin{lem}\label{1wg.17}  Let all the assumptions in Theorem \ref{wg.7.2_1} hold.   Let u solve the system (\ref{wg.1}).
 Then
\be \label{1wg.17.1}\int_0^T\int_{\Om}\left(u_t^2+u^2_\rho+(2m-1) |\nabla_{\Ga_g}u|^2_g\right) dx dt\leq \int_{\Om}\rho \left(u^2_t+ u^2_\rho\right)dx + C(R_0)E(0), \ee
where $\rho(x)$ is given by (\ref{cpde.2.2}).
\end{lem}

 {\bf Proof}.
Let $\rho(x)$ be given by (\ref{cpde.2.2}), then
\be  \rho(x)=r^m, x\in\overline{\Om}. \ee

 With (\ref{wg.7.2_3}) and (\ref{pde.5}), we have
\be D^2\rho(X,X) \ge \frac{m_1}{\rho}|X|_g^2,\quad \textmd{for}\ \  X\in S(r)_x, x\in\Om. \ee
Let $\H=\rho\frac{\partial}{\partial \rho}$ in (\ref{wg.14.1}). With (\ref{wg.7.2_3}) and (\ref{pde.5}),  we deduce that
\begin{eqnarray} \label{wg.16_1}D\H\left(\nabla_gu,\nabla_g u\right)=&&D\H\left(u_\rho\frac{\partial}{\partial \rho} ,u_\rho\frac{\partial}{\partial \rho}\right)+D\H(\nabla_{\Ga_g}u,\nabla_{\Ga_g}u)\nonumber\\
=&&u^2_\rho+\rho D^2\rho\left(\nabla_{\Ga_g}u,\nabla_{\Ga_g}u\right)=u^2_\rho+m_1 |\nabla_{\Ga_g}u|^2_g,\quad  x\in\Om,\end{eqnarray}
and
\be\label{wg.16_2}\quad \div \H=\div\rho \frac{\partial}{\partial \rho}=\frac{1}{m}\div r \frac{\partial}{\partial r}=\frac{n}{m}, \quad  x\in\Om.\ee

From  (\ref{wg.14.1}),  we obtain
\begin{eqnarray}
 \label{wg.16}
 &&\int_0^T\int_{\partial\Om(a)}\frac{\pa u}{\pa\nu_{\A}}\H(u) d\Ga dt+\frac12\int_0^T\int_{\partial\Om(a)}
\left(u_t^2-\left|\nabla_g u\right|_g^2\right)\H\cdot\nu d\Ga dt\nonumber\\
 =&&(u_t,\H(u))\Big |^T_0+\int_0^T\int_{\Om(a)}\left(u^2_\rho+m_1 |\nabla_{\Ga_g}u|^2_g\right) dx dt\nonumber\\
&&+\int_0^T\int_{\Om(a)}\frac{\div \H}{2}\left(u_t^2-\left |\nabla_g
u\right|_g^2\right) dx dt\nonumber\\
=&& (u_t,\H(u))\Big |^T_0+\frac{1}{2}\int_0^T\int_{\Om(a)}\left(u_t^2+u^2_\rho+(2m_1-1) |\nabla_{\Ga_g}u|^2_g\right) dx dt\nonumber\\
&&+\int_0^T\int_{\Om(a)}\frac{1}{2}\left(\frac{n}{m}-1\right)\left(u_t^2-\left |\nabla_g
u\right|_g^2\right) dx dt.
\end{eqnarray}
Let $P=\frac{1}{2}\left(\frac{n}{m}-1\right)$, substituting (\ref{wg.14.2}) into (\ref{wg.16}), letting $a\rightarrow\infty$,  we derive that
 \begin{eqnarray}
 \label{wg.16.1}
 &&\Pi_{\Ga}=\left(u_t,\H(u)+Pu\right)\Big|^T_0+\frac{1}{2}\int_0^T\int_{\Om(a)}\left(u_t^2+u^2_\rho+(2m_1-1) |\nabla_{\Ga_g}u|^2_g\right) dx dt,
\end{eqnarray}
where
\begin{eqnarray}
 \label{wg.16.1.111}
\Pi_{\Ga}=&&\int_0^T\int_{\Ga}\frac{\pa u}{\pa\nu_{\A}}\H(u) d\Ga dt+\frac12\int_0^T\int_{\Ga}
\left(u_t^2-\left|\nabla_g u\right|_g^2\right)\H\cdot\nu d\Ga dt\nonumber\\
&&+\int_0^T\int_{\Ga}Pu\frac{\pa
u}{\pa\nu_{\A}}d\Ga dt-\frac12\int_0^T\int_{\Ga}u^2\frac{\pa
P}{\pa\nu_{\A}}d\Ga dt.
\end{eqnarray}

 Since $u\Large|_{\Ga}=0,$ we
obtain $\nabla_{\Ga_g} u\large|_{\Ga}=0$, that is,
 \be\nabla_g u=\frac{\pa
 u}{\pa\nu_{\A}}\frac{\nu_{\A}}{|\nu_{\A}|_g^2}\quad\mbox{for}\quad
 x\in\Ga.\label{wg.16.1.2}\ee
Similarly, we have \be \H(u)=\<\H,\nabla_g u\>_g=\frac{\pa
 u}{\pa\nu_{\A}}\frac{\H\cdot\nu}{|\nu_{\A}|_g^2}\quad\mbox{for}\quad
 x\in\Ga.\label{wg.16.1.3}\ee
With (\ref{wg.7_1}) and  (\ref{wg.7_2}), we deduce that
\be\label{wg.16.1_2} \left\<\frac{\partial}{\partial \rho},\nu\right\>=\left\<\frac{\partial}{\partial r}\frac{dr}{ d \rho},\nu\right\>=\frac{r^{1-m}}{m} \left\<\frac{\partial}{\partial r},\nu\right\>=\frac{r^{1-m}}{m} \frac{\partial r}{\partial \nu}\le 0 \quad \mbox{for} \ \ x\in \Ga.\ee
Using the formula (\ref{wg.16.1.3}) in the formula
(\ref{wg.16.1.111}) on the portion $\Ga$, with (\ref{wg.16.1_2}), we obtain
\begin{equation}
\label{wg.16.1.4} \Pi_{\Ga}=\frac12\int_0^T
\int_{\Ga}\left(\frac{\pa
u}{\pa\nu_{\A}}\right)^2\frac{\H\cdot\nu}{|\nu_{\A}|_g^2} d\Ga dt\le 0.
\end{equation}

Substituting (\ref{wg.16.1.4}) into (\ref{wg.16.1}),we have
 \begin{eqnarray}
 \label{wg.16.2}
 &&\left(u_t,\rho u_\rho+\frac{1}{2}\left(\frac{n}{m}-1\right) u\right)\Big|^T_0\nonumber\\
&& \qquad+\frac{1}{2}\int_0^T\int_{\Om}\left(u_t^2+u^2_\rho+(2m_1-1) |\nabla_{\Ga_g}u|^2_g\right) dx dt\le0.
\end{eqnarray}
With  (\ref{wg.17.1}) and (\ref{wg.16.2}), we deduce that
\begin{eqnarray}
 \label{wg.20.5}
 &&\int_0^T\int_{\Om}\left(u_t^2+u^2_\rho+(2m_1-1) |\nabla_{\Ga_g}u|^2_g\right) dx dt\nonumber\\
 \leq&&  2\int_{\Om} \rho \left|u_t\left( u_\rho +\frac{1}{2\rho}\left(\frac{n}{m}-1\right)u\right)\right| dx + C(R_0)E(0)\nonumber\\
 \leq&& \int_{\Om} \rho \left(u^2_t+ u^2_\rho +\frac{1}{\rho}\left(\frac{n}{m}-1\right)uu_\rho+\frac{1}{4\rho^2}\left(\frac{n}{m}-1\right)^2u^2\right)dx + C(R_0)E(0)\nonumber\\
 =&& \int_{\Om} \rho\left(u^2_t+ u^2_\rho\right)dx +\int_{\Om} \left(\frac{1}{\rho}\left(\frac{n}{m}-1\right)uu_\rho+\frac{1}{4\rho^2}\left(\frac{n}{m}-1\right)^2u^2\right)dx+ C(R_0)E(0) \nonumber\\
=&&  \int_{\Om}  \rho\left(u^2_t+ u^2_\rho\right)dx -\int_{\Om} \frac{1}{4\rho}\left(\frac{n}{m}-1\right)^2u^2 dx+ C(R_0)E(0) \nonumber\\
\leq&& \int_{\Om}  \rho\left(u^2_t+ u^2_\rho\right)dx+ C(R_0)E(0) .
\end{eqnarray}
Then, the estimate (\ref{1wg.17.1}) holds.
$\Box$
$$ $$

{\bf Proof of Theorem \ref{wg.7.2_1} }\,\,\,


 Substituting (\ref{1wg.20}) into  (\ref{1wg.17.1}), we have
\be \label{1wg.20.1}\int_0^T\int_{\Om} |\nabla_{\Ga_g}u|^2_g dx dt\leq   C(R_0)E(0).  \ee
 Substituting (\ref{1wg.20.1}) into  (\ref{1wg.17.1}), we obtain
\be TE(T)\leq \int_{\Om}\rho(u_t^2+u_\rho^2)dx + C(R_0)E(0)\leq  \int_{\Om \backslash\Om(a) }\rho\left(u_t^2+\left |\nabla_g
u\right|_g^2\right)dx +C(a,R_0)E(0).  \ee
With (\ref{1wg.20.3}), we deduce that
\begin{eqnarray}TE(T,a)&&\leq  \int_{\Om \backslash\Om(a) }(\rho-T)\left(u_t^2+\left |\nabla_g
u\right|_g^2\right)dx +C(a,R_0)E(0)\nonumber\\
&& \leq \int_{\Om  }e^{\rho-T}\left(u_t^2+\left |\nabla_g
u\right|_g^2\right)dx +C(a,R_0)E(0)\leq C(a,R_0)E(0).  \end{eqnarray}

The estimate (\ref{wg.7.2_2}) holds.

\end{document}